\documentclass[a4paper]{jpconf}
\usepackage{graphicx}
\usepackage{blopeps}

\usepackage[thmmarks]{ntheorem-hyper}
\theoremheaderfont{\bf}\theorembodyfont{\it}\theoremsymbol{}\theoremstyle{plain}\newtheorem{Thm}{Theorem}
\theoremheaderfont{\bf}\theorembodyfont{\rm}\theoremsymbol{}\theoremstyle{plain}\newtheorem{Prf}{Proof}
\theoremheaderfont{\bf}\theorembodyfont{\it}\theoremsymbol{}\theoremstyle{plain}
\theoremheaderfont{\bf}\theorembodyfont{\sl}\theoremsymbol{}\theoremstyle{plain}
\theoremheaderfont{\bf}\theorembodyfont{\it}\theoremsymbol{}\theoremstyle{plain}
\theoremheaderfont{\bf}\theorembodyfont{\sl}\theoremsymbol{}\theoremstyle{plain}
\theoremheaderfont{\bf}\theorembodyfont{\it}\theoremsymbol{}\theoremstyle{plain}

\begin{document}
\title{A linear iterative unfolding method}

\author{Andr\'as L\'aszl\'o}

\address{Wigner RCP, P.O.Box 49, H-1525 Budapest, Hungary\newline CERN, CH-1211 Geneve 23, Switzerland}

\ead{laszlo.andras@wigner.mta.hu}

\begin{abstract}
A frequently faced task in experimental physics is to measure the 
probability distribution of some quantity. Often this quantity to be 
measured is smeared by a non-ideal detector response or by some physical 
process. The procedure of removing this smearing effect from the measured 
distribution is called unfolding, and is a delicate problem in signal 
processing, due to the well-known numerical ill behavior of this task. 
Various methods were invented which, given some assumptions on the initial 
probability distribution, try to regularize the unfolding problem. 
Most of these methods definitely introduce bias into the estimate of the 
initial probability distribution. We propose a linear iterative method 
(motivated by the Neumann series / Landweber iteration known in functional 
analysis), which has the advantage that no assumptions on the initial 
probability distribution is needed, and the only regularization parameter 
is the stopping order of the iteration, which can be used to choose the best 
compromise between the introduced bias and the propagated statistical and 
systematic errors. The method is consistent: 
``binwise'' convergence to the initial probability 
distribution is proved in absence of measurement errors under a quite general 
condition on the response function. This condition holds for practical applications 
such as convolutions, calorimeter response functions, momentum reconstruction 
response functions based on tracking in magnetic field etc. In presence of measurement 
errors, explicit formulae for the propagation of the three important 
error terms is provided: bias error (distance from 
the unknown to-be-reconstructed initial distribution at a finite iteration 
order), statistical error, and systematic error. 
A trade-off between these three error terms can be used to define an optimal 
iteration stopping criterion, and the errors can be estimated there. We provide a 
numerical C library for the implementation of the method, which incorporates 
automatic statistical error propagation as well. 
The proposed method is also discussed in the context of other known approaches.
\end{abstract}

\section{Introduction}

In data analysis one commonly faces the problem that the probability density 
function (pdf) of a given physical quantity of interest is to be measured, 
but some random physical process, such as the intrinsic behavior of the 
measurement apparatus, smears it. The reconstruction of the pertinent pdf 
based on the measured smeared pdf and on the response function of the 
measurement procedure is called unfolding. To be specific, let us have 
the original unknown pdf $x\mapsto f(x)$ of the undistorted physical quantity 
which we need to reconstruct, 
and assume that the actual measured pdf can be expressed of the form 
$y\mapsto g(y)=\int\rho(y|x)\,f(x)\,\mathrm{d}x$, where $(y,x)\mapsto\rho(y|x)$ 
describes the smearing effect in a probabilistic manner.\footnote{All pdfs are understood 
to be real valued non-negative Lebesgue integrable functions over some 
finite dimensional real vector space $X$.} 
Then, it is said that the pdf $g$ is the pdf $f$ folded with 
the response function $\rho$.\footnote{
Whenever the response function $\rho$ is translation invariant in the 
sense that for all $x,y,z\in X$ one has $\rho(y+z|x)=\rho(y|x-z)$, the folding 
is specially called convolution, and in that case $\rho$ may be expressed 
by a single pdf: $\rho(y|x)=\rho(y-x|0)$.} 
Our mathematical task is to solve the above linear integral equation in 
order to obtain $f$, given $g$ and $\rho$. This problem is known not to 
be a simple numerical task (ill-posed problem), and several methods are used 
by the data analysis communities in order to regularize the problem (for 
an overview on the most popular approaches, we refer to \cite{cowan2002, blobel2008}).

\begin{figure}[!ht]
\begin{center}
\includegraphics[width=12cm]{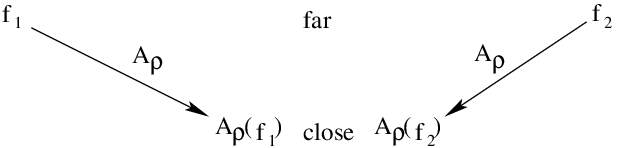}
\end{center}
\caption{\label{closedgraph} Illustration of the non-continuity of the 
inverse of a folding operator $A_{\rho}$: two distant functions $f_{1}$ and 
$f_{2}$ may be mapped close by the folding -- distance of functions are here 
understood as probabilistic distance, i.e.\ in the $L^{1}(X)$ function norm.}
\end{figure}

Let us denote by $A_{\rho}$ the pertinent folding operator, which acts like 
$\left(A_{\rho}f\right)(y)=\int\rho(y|x)\,f(x)\,\mathrm{d}x$ on a function 
$f$ at a point $y$.\footnote{To be precise, 
$A_{\rho}$ is a $L^{1}(X)\rightarrow L^{1}(X)$ continuous linear operator, 
where $L^{1}(X)$ denotes the normed space of complex valued integrable 
functions over the vector space $X$. The response function $\rho$ is assumed 
to be $\rho(\cdot|x)\in L^{1}(X)$ for all $x\in X$.}
Given the measured pdf $g=A_{\rho}f$, the problem of unfolding 
can then be formalized as follows: the pdf $f=A_{\rho}^{-1}(g)$ is to 
be determined or approximated. The mathematical cause of the numerical 
ill-posedness of this unfolding problem can then be put forward as: 
the inverse $A_{\rho}^{-1}$ of a generic folding operator can be shown not 
to be continuous despite the forward folding operator $A_{\rho}$ always 
being continuous\footnote{Continuous in the 
$L^{1}(X)\rightarrow L^{1}(X)$ sense.} 
(this phenomenon is discussed in detail e.g.\ in \cite{laszlo2006}). 
The non-continuity of the inverse folding operator $A_{\rho}^{-1}$ may be 
also reformulated in a less abstract manner: initially distant functions 
can be mapped close by the folding operator $A_{\rho}$, as 
illustrated in Figure~\ref{closedgraph}. I.e.\ one can lose discriminating 
power between pdfs upon a folding.

\begin{figure}[!ht]
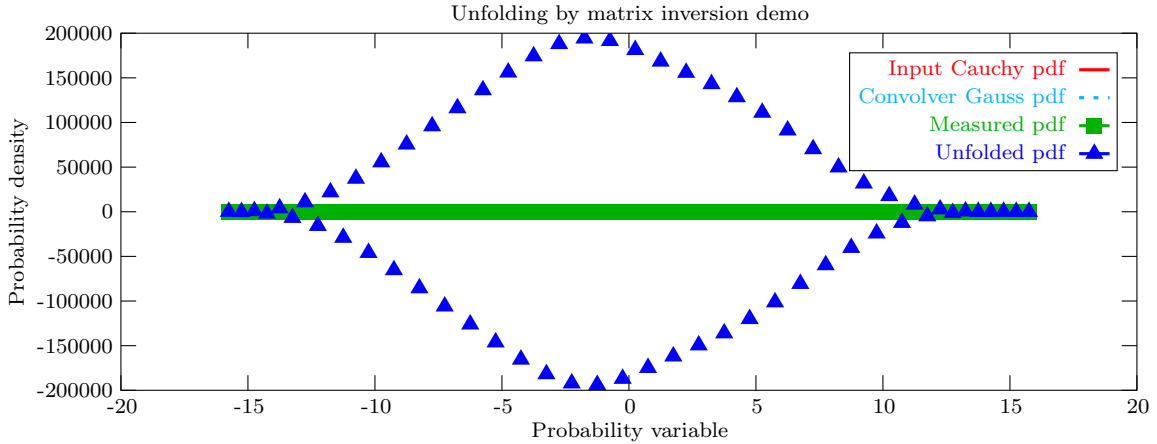

\begin{center}
{\scriptsize\blopeps[width=15cm,height=6cm]{fig/inverse_gausscauchy/unfoldedpdf_inverse.beps}}
\end{center}
\caption{\label{nononto} (Color online) Demonstration of the numerical ill-posedness 
of the unfolding problem: a Cauchy distribution is convolved with a Gauss 
distribution with Monte Carlo method to generate the measured distribution 
contaminated with statistical errors. Clearly, the unfolded pdf, 
obtained by simple numerical inversion of the discretized 
folding operator on the measured pdf gives physically unreasonable numerical 
result: large alternating positive / negative amplitude pdf values.}
\end{figure}

A further aspect of the numerical ill-posedness of the unfolding problem is that 
in practice the folded pdf $g$ is often obtained via statistical measurements 
(e.g.\ histograming), and therefore is contaminated by statistical errors. 
I.e.\ in reality $g=A_{\rho}f+e$ holds instead of the idealized equation 
$g=A_{\rho}f$, where $e(x)$ is a random variable for each point 
$x$ (or for each histogram bin -- in the language of histograms). Thus, 
when estimating the unfolded pdf as $A_{\rho}^{-1}(g)=f+A_{\rho}^{-1}(e)$, 
the contribution of the second term is not guaranteed to remain small due to 
the non-continuity of the inverse folding operator $A_{\rho}^{-1}$ 
even when $e$ is initially known to be small. On top of this, the statistical 
error term $e$ may contain modes not within the image of the folding operator 
$A_{\rho}$, on which the evaluation of the inverse operator $A_{\rho}^{-1}$ 
is not meaningful if the problem is not initially discrete. These effects are demonstrated 
in Figure~\ref{nononto}, which shows that simple inversion of the discretized 
folding operator on the measured pdf gives unphysical numerical result: a result 
very different from the initial pdf, having large negative and positive alternating amplitudes.

In order to regularize the numerical ill-posedness of the unfolding problem, 
various methods are used. These methods can be divided into three large classes.
\begin{enumerate}
 \item Using a parametric Ansatz for $f$, and fit parameters, so that $A_{\rho}f$ 
 gets close to $g$. This method can be slightly insensitive to the details of the 
 true $f$ (as illustrated in Figure~\ref{closedgraph}), and of course can introduce 
 strong systematic bias on the result if the parametric Ansatz does not hold 
 in an exact manner of the form that was assumed. Such methods are used in general 
 for inclusive particle identification by specific ionization (see e.g.\ \cite{laszlo2008}).
 \item Bin-by-bin fitting of the bin values of the histogramed $f$, so that $A_{\rho}f$ 
 gets close to $g$. This is 
 basically equivalent to the naive inversion of the discretized folding operator, 
 and therefore produces similar oscillatory results, except when an artificial 
 penalty function is added to the $\chi^{2}$ in order to suppress large local gradients. In that 
 case, the method can provide meaningful answers, but the introduced systematic 
 bias is difficult to quantify. Similarly to the parametric Ansatz method, 
 the fit can be slightly insensitive to the details of the true $f$. 
 Most popular methods, such as SVD method \cite{hoecker1996}, are based on this 
 idea.
 \item The iterative method of convergent weights 
 (also known as iterative Bayesian unfolding) of 
 Kondor-M\"ulthei-Schorr-d'Agostini 
 \cite{shepp1982, kondor1983, multhei1987_1, multhei1987_2, dagostini1995}. 
 This method is, as opposed to the previously mentioned methods, is non-linear. 
 On the other hand, by construction it preserves positivity and 
 integral of the initial pdf, and therefore maps a pdf exactly into a pdf, 
 which does not hold for linear methods, thus, this approach is quite 
 favorable for statistical applications. Regularization is achieved solely by 
 stopping the iteration at a finite order. However, there is no known proof yet 
 if the iterated pdfs converge\footnote{In case of an iterative unfolding method 
 it is an absolute must to show that the sequence of iterated unfolded pdfs 
 converge to the initial one, in the absence of measurement errors (consistency of the method).} 
 to the initial to-be-reconstructed pdf in a non-discrete 
 scenario, even in the absence of measurement errors \cite{multhei1987_1}. Also 
 propagation of statistical and systematic errors of the measured pdf to the 
 unfolded pdf has not been investigated, and consequently 
 no generally applicable iteration stopping condition is known.
\end{enumerate}
In a previous paper \cite{laszlo2006} we proposed a linear iterative unfolding 
method, for which under certain conditions convergence to the initial pdf 
was proved analytically for some unfolding problems in probability 
theory (such as convolutions), and due to the linearity of the method, 
exact propagation of statistical errors of the measured 
(folded) pdf to the unfolded pdf was possible. 
In this paper we propose an improved version of that algorithm, which could be 
proved to be convergent in quite general cases for unfolding problems in a 
probability theory setting.\footnote{The detailed mathematical proof of convergence 
shall be published elsewhere: \cite{laszloinpr}. The proposed iteration scheme was motivated by the so 
called Neumann series and Landweber iteration \cite{landweber1951} 
known in functional analysis, but the convergence of neither iterations 
hold, unfortunately, in a probability theory setting in their original form, 
as one can prove. Our improved iterative algorithm, however, 
is specially developed to be convergent for unfolding problems in probability 
theory.} The key equality of the 
convergence proof leads to explicit error propagation formulae for the three 
important error terms: for the bias error (distance from the true unfolded pdf), 
for the propagated statistical error, and most notably for the propagated 
systematic error, which is of great importance in reporting experimental results. 
An implementation of the algorithm is written as C library, along with 
application examples \cite{laszlo2011}. The implementation also incorporates 
automatic statistical error propagation.

The paper is organized as follows: in Section~\ref{convergence} the algorithm 
and its convergence theorem shall be formulated, Section~\ref{errorpropagation} 
is devoted to the corresponding error propagation formulae which help to 
formulate an optimal stopping criterion and error estimates therein, while 
in Section~\ref{examples} we demonstrate our method on examples.

\section{A linear iterative unfolding algorithm}
\label{convergence}

We provide now a linear iterative solution for a probability theory unfolding 
problem of the form $g=A_{\rho}f$, where $f$ is the initial (unknown) pdf, 
$g$ is the folded (measured) pdf, and $\rho$ is the response function. 
Given the response function $\rho$, one can also define along with the folding 
operator $A_{\rho}$ the transpose folding operator $A_{\rho}^{T}$ by swapping 
the variables of the response function.\footnote{The transpose folding 
operator is defined by $(A_{\rho}^{T}f)(y)=\int\rho(x|y)\,f(x)\,\mathrm{d}x$ 
for all functions $f$ and points $y$. Note, that this simply translates to 
matrix transposition whenever the folding is discretized.} Then, one can attempt 
to approximate the true unfolded pdf $f$ in the following way:\newline
define the function sequence by setting the normalization factor
\begin{equation}
K_{\rho}=\max_{x}\,\int\int \,\rho(y|z) \,\rho(y|x) \,\mathrm{d}y\,\mathrm{d}z
\end{equation}
and then taking the 
\begin{eqnarray}
f_{0}  &=&K_{\rho}^{-1}A_{\rho}^{T}g,\cr
f_{N+1}&=&f_{N}+\left(f_{0}-K_{\rho}^{-1}A_{\rho}^{T}A_{\rho}f_{N}\right)
\end{eqnarray}
iteration formula. We provide a convergence result on this iterative approximation 
below in absence of measurement errors on $g$ (which is necessary for the consistency of the method).

\begin{figure}[!ht]
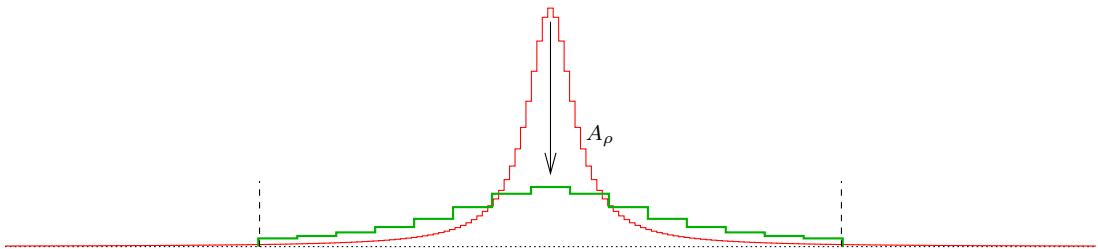

\begin{center}
{\scriptsize\blopeps[width=15cm,height=4cm]{fig/rebinning/rebinning.beps}}
\end{center}
\caption{\label{rebinning} (Color online) Illustration the rebinning trick 
for unfolding the histogram binning and domain truncation as well 
(to the maximum possible extent) along with the smearing effect of the 
response function $\rho$. 
For this, the implementation of the folding operator $A_{\rho}$ must map histograms 
over a larger domain and with finer graining to histograms with the binning 
scheme of the measured (folded) pdf.}
\end{figure}

\begin{Thm}
\label{unfoldingtheorem}
(Convergence) 
The function sequence $N\mapsto f_{N}$ resulting from the above iteration scheme 
converges to the closest possible function to the true unfolded pdf $f$ in the 
average over any compact region, whenever the normalization factor $K_{\rho}$ 
is finite. I.e.\ for all compact sets $S\subset X$ one has
\begin{equation}
\lim_{N\rightarrow\infty} \frac{1}{\mathrm{Volume}(S)}\int_{S} \left( f-P_{\mathrm{Ker}(A_{\rho})}f - f_{N} \right)(x)\,\mathrm{d}x=0.
\end{equation}
Here, $P_{\mathrm{Ker}(A_{\rho})}$ denotes the orthogonal projection operator to 
the kernel set of $A_{\rho}$, and thus $P_{\mathrm{Ker}(A_{\rho})}=0$ holds 
automatically whenever $A_{\rho}$ is invertible. In 
addition, the convergence shall also hold in the space of square-integrable functions, 
i.e.\ one has also
\begin{equation}
\lim_{N\rightarrow\infty} \int \left\vert f-P_{\mathrm{Ker}(A_{\rho})}f - f_{N} \right\vert^{2}(x)\,\mathrm{d}x=0.
\end{equation}
\end{Thm}
\begin{Prf}
The proof is based on Riesz-Thorin theorem and on the spectral representation 
of positive operators in the space of complex square-integrable functions over 
$X$ (to be published in a more mathematically specialized journal: \cite{laszloinpr}).
\end{Prf}
The following observations help to shed some further light on properties of 
the proposed unfolding algorithm.
\begin{enumerate}
 \item In case the pdfs are modeled with histograming, the setwise convergence 
 of pdfs means binwise convergence of histograms, i.e. the probability of each 
 histogram bin is restored in the limit of infinite iterations.
 \item When the inverse of $A_{\rho}$ exists, the original pdf $f$ is completely 
 restored. Whenever the pertinent inverse does not exist, still the maximum 
 possible information about $f$ is restored, namely the function 
 $f-P_{\mathrm{Ker}(A_{\rho})}f$.
 \item Whenever $A_{\rho}$ is a convolution, then $K_{\rho}=1$ holds automatically, i.e.\ $K_{\rho}<\infty$ is satisfied.
 \item The convergence condition $K_{\rho}<\infty$ holds provably for 
 a wide class of practically relevant response functions, such as 
 energy response function of calorimeters, 
 momentum response function of track reconstruction in magnetic field etc.
 \item The iteration scheme of the theorem is motivated by the Neumann series 
 known in functional analysis. 
 A similar iterative solution, also referred to 
 as Landweber iteration \cite{landweber1951}, is known in the theory of 
 Fredholm operators. In probability theory unfolding problems, however, 
 the necessary convergence criteria for Neumann series or for Landweber 
 iteration do not hold in their original form.
 \item The proposed iterative unfolding algorithm does not necessarily need 
 an initial binning of pdfs. It may be implemented as well by different 
 density estimators than histograms. However, when the pdfs are modeled 
 by histograms, one may recognize that the binning and truncation of 
 histograming domain can also be considered  as folding operator. 
 Therefore, the histogram binning and truncation effect 
 may be included in the response function $\rho$, and then the effect of 
 histograming can be unfolded (to the maximum possible extent) as well. If 
 one wants to numerically implement this, the initial pdf $f$ must be assumed 
 to better approach the continuum pdf, i.e. must be assumed to be an unknown 
 histogram over a larger domain with finer granulation than the folded one. 
 The schematic of such possible rebinning trick is illustrated in Figure~\ref{rebinning}.
\end{enumerate}

\section{Bias, statistical and systematic errors of the unfolded distribution}
\label{errorpropagation}

In real measurements, the folded pdf $g$ also admits statistical and systematic 
errors, and the propagation of these terms into the unfolded pdf is necessary 
to quantify at each finite iteration step. The key equality of the proof of Theorem~1 leads to 
explicit error propagation formulae for bias error (distance from the true unfolded pdf), 
statistical error, and systematic error. First we present our result about 
bias error.

\begin{Thm}\label{biasthm} (Bias error)
Take the iterative solution for the unfolding problem as in Section~\ref{convergence}. Then, 
if the normalization factor $K_{\rho}$ is finite, the distance of an $N$-th 
iterate $f_{N}$ from the closest possible function to the true unfolded pdf 
$f$ in the average over a compact region has the following upper bound:
for any compact set $S\subset X$ one has
\begin{eqnarray}
\left\vert\frac{1}{\mathrm{Volume}(S)}\int_{S}\left(f-P_{\mathrm{Ker}(A_{\rho})}f-f_{N}\right)(x)\,\mathrm{d}x\right\vert \leq \frac{1}{\sqrt{\mathrm{Volume}(S)}}\,(1+\varepsilon)\,\sqrt{\int \left\vert f_{M}-f_{N} \right\vert^{2}(x) \,\mathrm{d}x}\cr
\end{eqnarray}
for any $\varepsilon>0$ and large enough iteration order $M>N$.
\end{Thm}

The above result, translated to the language of histograms means that the 
bin-by-bin average deviation from the true unfolded pdf is bound by the 
right hand side of the inequality in Theorem~\ref{biasthm}, 
where $\mathrm{Volume}(S)$ is the histogram bin volume, $N$ is the iteration order, and 
$(1+\varepsilon)\,\sqrt{\int \left\vert f_{M}-f_{N} \right\vert^{2}(x) \,\mathrm{d}x}$ 
is a calculable coefficient. In this expression $\varepsilon>0$ is arbitrary, 
while the iteration order $M>N$ needs to be large enough for given $\varepsilon$. 
It is seen that the bias error tends to zero with increasing iteration order $N$ 
and depends on the histogram bin size as $\frac{1}{\sqrt{\mathrm{Volume}(S)}}$.

In practical applications, the pdfs are often measured by statistical methods (e.g.\ histograming). 
In that case, the value of the folded pdf $g$ in each histogram bin admits a 
statistical error. The below theorem states an exact formula for the propagation 
of this error into the unfolded pdf.

\begin{Thm} (Statistical error)
Take the iterative solution for the unfolding problem as in Section~\ref{convergence}. 
Let $C$ be the covariance matrix of the measured pdf $g$, where $g$ is assumed to be of the 
form of a histogram. Since a covariance matrix $C$ is positive definite, 
it is always possible to decompose it -- not uniquely -- in the form 
$C=E\,E^{T}$ for some matrix $E$, $(\cdot)^{T}$ being the matrix transpose. 
(Whenever $C$ is diagonal, construction of such an $E$ is just trivial.)
Then, the following iteration calculates the statistical error propagation:
\begin{eqnarray}
E_{0}  &=&K_{\rho}^{-1}A_{\rho}^{T}E,\cr
E_{N+1}&=&E_{N}+\left(E_{0}-K_{\rho}^{-1}A_{\rho}^{T}A_{\rho}E_{N}\right),
\end{eqnarray}
where in each step the covariance matrix of $f_{N}$ shall be $C_{N}=E_{N}\,E_{N}^{T}$.
\end{Thm}

Due to the linearity of the method, the contribution of the propagated 
statistical error term is exactly calculable by means of the above formulae, 
if it is known for the measured pdf $g$. This error term increases with 
increasing iteration order $N$. The statistical error of a given histogram bin 
of the $N$-th iterate $f_{N}$ is nothing but the square-root of the corresponding 
diagonal element of $C_{N}$.

Whenever the folded pdf $g$ is a result of an experiment, it may admit a systematic 
error $\delta{g}$. Also the systematic error $\delta{\rho}$ of the response 
function $\rho$ may give a non-zero contribution to it: $A_{\delta{\rho}}f$. 
The effect of this initial systematic error on the unfolded pdf is quantified 
by the following theorem.

\begin{Thm}\label{systthm} (Systematic error)
Take the iterative solution for the unfolding problem as in Section~\ref{convergence}. 
Assume that $\delta{g}$ is the systematic error of $g$ (possibly including 
contribution from systematic error of the response function). Then the 
systematic error for the $N$-the iterate $f_{N}$ averaged over a compact 
region has the following upper bound: for any compact set $S\subset X$
\begin{eqnarray}
\left\vert\frac{1}{\mathrm{Volume}(S)}\int_{S} \delta{f}_{N}(x)\,\mathrm{d}x\right\vert \leq \sqrt{\int \left\vert \Xi_{{}_{S,N}}\right\vert^{2}(x) \,\mathrm{d}x}\, \sqrt{\int \left\vert K_{\rho}^{-1}A_{\rho}^{T}\delta{g}\right\vert^{2}(x)\,\mathrm{d}x}\cr
\end{eqnarray}
where $\Xi_{{}_{S,N}}$ is defined by the iteration
\begin{eqnarray*}
\Xi_{{}_{S,0}}   & := & \frac{1}{\mathrm{Volume}(S)} \chi_{{}_{S}}, \cr
\Xi_{{}_{S,N+1}} & := & \Xi_{{}_{S,N}} + \left(\Xi_{{}_{S,0}} - K_{\rho}^{-1}A_{\rho}^{T}A_{\rho}\Xi_{{}_{S,N}}\right)\cr
\end{eqnarray*}
$\chi_{{}_{S}}$ being the characteristic function of the set $S$.
\end{Thm}

The above result, translated to the language of histograms means that the 
bin-by-bin average systematic error of the $N$-th iterate $f_{N}$ is bound by the 
formula in the right hand side of the inequality in Theorem~\ref{systthm}, 
where $\mathrm{Volume}(S)$ is the histogram bin volume, $N$ is the iteration order, and the 
coefficient $\sqrt{\int \left\vert K_{\rho}^{-1}A_{\rho}^{T} \delta{g} \right\vert^{2}(x) \,\mathrm{d}x}$ 
is calculable knowing the bin-by-bin systematic errors $\delta{g}$ of the measured pdf $g$. 

As the bias error decreases, while the statistical and systematic error of 
the $N$-th iterate $f_{N}$ increases with the iteration order $N$, a trade-off 
between these error terms provides an optimal cutoff criterion\footnote{E.g.\ 
one can take the sum of the three error terms, and stop the iteration when it 
reaches a minimum.} in the iteration order $N$, and error estimates therein. 
Consequently the true unfolded pdf $f$ can be approximated optimally and 
the error of this approximation can be put under full control. 
Thus, the regularization of the numerically ill-posed unfolding problem is 
achieved, in case of the proposed approach, solely by using an iterative 
approximation and choosing an optimal iteration stop order taking into account 
the convergent terms (bias error) and the divergent terms (statistical and 
systematic errors).

\section{Examples}
\label{examples}

\begin{figure}[!ht]
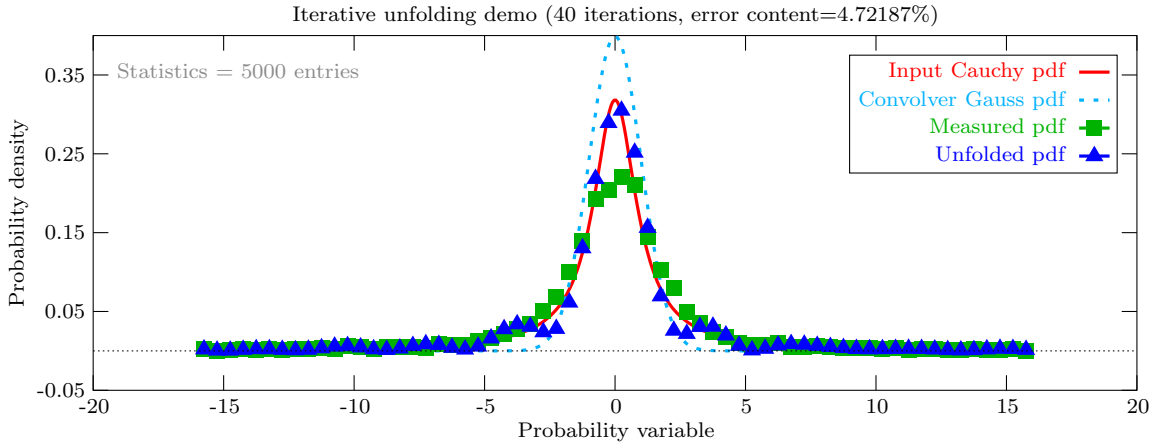

\begin{center}
{\scriptsize\blopeps[width=15cm,height=6cm]{fig/iterative_gausscauchy/unfoldedpdf_iterative.beps}}
\end{center}
\caption{\label{gausscauchy} (Color online) Test example with unfolding a Cauchy distribution 
convolved with a Gaussian distribution. Iteration was stopped when the integral 
of the statistical error term reached about $5\%$.}
\end{figure}

In this section we give two examples to demonstrate our method. In the 
first example, we take a Cauchy distribution, and convolve it with a 
Gaussian distribution with Monte Carlo method. The folded pdf is determined 
by histograming the sum of the Cauchy and Gaussian distribution random numbers, 
i.e.\ the measured pdf shall admit Poissonian statistical errors. In this example, a 
relatively modest statistics of 5000 entries was taken to be able to judge the 
method in the low statistics limit. The results is shown 
in Figure~\ref{gausscauchy}. It is seen that the original Cauchy pdf is restored, 
modulo the fluctuations arising from the propagated statistical errors -- these 
are seen as ``shoulders'' of the unfolded pdf, the amplitude of which decrease 
with increased statistics. The iteration was stopped when the integral of 
the statistical error term reached about $5\%$ level.

\begin{figure}[!ht]
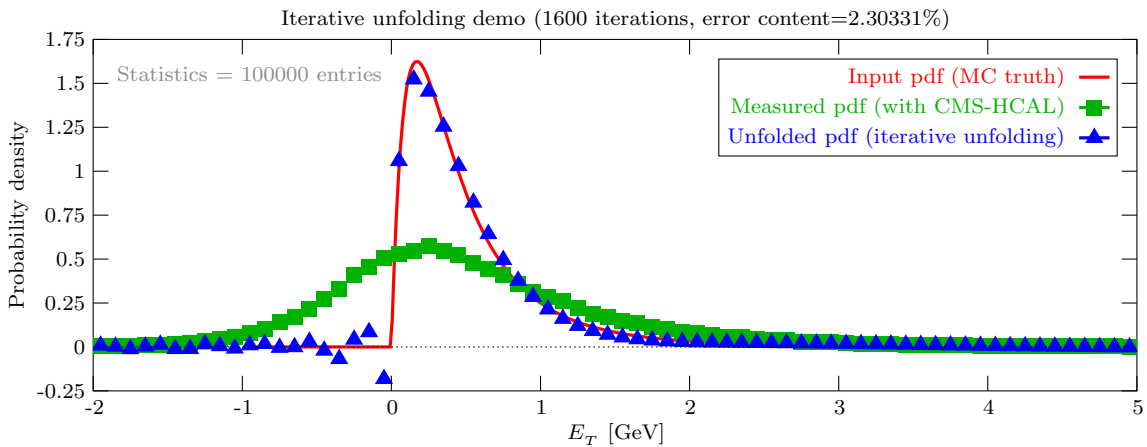

\begin{center}
{\scriptsize\blopeps[width=15cm,height=6cm]{fig/iterative_cmshcaltsallis/cmshcaltsallis_iterative.beps}}
\end{center}
\caption{\label{cmshcal} (Color online) A physical example with unfolding 
energy distribution of charged hadrons measured with hadronic calorimeter.
Iteration was stopped when the integral of the statistical error term reached 
about $2.3\%$.}
\end{figure}

In the second Monte Carlo example, we generate the energy distribution of transversely emitted 
hadrons in $7\,\mathrm{GeV}$ p+p collisions \cite{khachatryan2010}, and 
we assume that this particle spectrum was measured by the CMS-HCAL calorimeter 
\cite{yazgan2009}. The unfolded spectrum, along with the true and measured 
distribution is shown in Figure~\ref{cmshcal}.

\section{Concluding remarks}

We proposed a linear iterative spectrum unfolding method for application 
in data analysis. Convergence to the true unfolded pdf is proved under a quite 
general condition \cite{laszloinpr} in absence of measurement errors, and 
error propagation formulae are derived for bias error, statistical 
error, and systematic error in the presence of measurement errors. The method 
is demonstrated on physical examples. A numerical library in C is provided 
with the implementation of the method \cite{laszlo2011}. The algorithm could 
be included in the ROOUnfold package \cite{adye2011} in the future.

\ack

I would like to thank prof. G\"unter Zech for valuable discussions.

\end{document}